\documentclass[11pt]{article}

\usepackage{amsmath,amssymb,xfrac}
\usepackage{algorithm,algorithmic}
\usepackage[margin=1in]{geometry}
\usepackage{url}
\usepackage{tikz}
\usetikzlibrary{arrows,calc,positioning,shapes.misc}

\usepackage{amsthm}

\newtheorem{assumption}{Assumption}
\theoremstyle{definition}

\usepackage{CustomCommands}

\title{A note on generalized semi-infinite program bounding methods}
\author{Stuart M. Harwood%
\thanks{%
Corporate Strategic Research \newline
ExxonMobil Research and Engineering \newline
Annandale, NJ 08801 USA \newline
{stuart.m.harwood@exxonmobil.com}\newline
ORCID: 0000-0001-5883-9624
}}
\date{\today}

\begin{document}
\maketitle

\begin{abstract}
Generalized semi-infinite programs (GSIP) are a class of mathematical optimization problems that generalize semi-infinite programs, which have a finite number of decision variables and infinite constraints.
Mitsos et al.~\cite{mitsosEA15} present a method for global optimization of GSIP.
This method involves a lower bounding method, and they claim that these lower bounds converge to the optimal objective value of the GSIP.
A counterexample is presented that shows that this claim is false.
\end{abstract}

\section{Introduction}
This note discusses methods for the global solution of generalized semi-infinite programs (GSIP).
Specifically, the method from \cite{mitsosEA15} is considered, and it is shown with counterexamples that the lower bounds do not always converge. 

Consider a GSIP in the general form
\begin{alignat}{2} 
\tag{GSIP}
\label{eq:gsip}
f^* = 
\inf_{x}\; & f(x) \\
\st 				
\notag & x \in X, \\
\notag & 0 \le \inf\set{ g(x,y) : y \in Y , h_j(x,y) \le 0, \forall j \in J },
\end{alignat}
for subsets $X$, $Y$ of finite dimensional real vector spaces, real-valued functions $f$, $g$, and $h_j$, and some finite index set $J$.
The last constraint (the infinite constraint) of \eqref{eq:gsip}
can be written in a number of ways, highlighting different aspects of the problem:
\begin{align}
\notag
&0 \le g(x,y), \;\forall y \in Y : h_j(x,y) \le 0, \forall j;\\
\notag
&\big([0 \le g(x,y)] \lor [\exists j : h_j(x,y) \nleq 0]\big), \;\forall y \in Y; \\
\notag
&\big([0 \le g(x,y)] \lor [\max_j\{ h_{j}(x,y) \} > 0]\big), \;\forall y \in Y.
\end{align}

The popular approach taken in the literature is to look at the problem%
\footnote{In many cases we may drop outer brackets grouping clauses in a disjunction.
Consider Boolean-valued $A$ and $B$;
if for all $y \in Y$ either $A(y)$ or $B(y)$ holds, we write
\[
A(y) \lor B(y), \;\forall y \in Y.
\]
}
\begin{alignat}{2} 
\label{eq:gsip_rel}
f^L = 
\inf_{x}\; & f(x) \\
\st 				
\notag & x \in X, \\
\notag & [g(x,y) \ge 0] \lor [\max_j\{ h_{j}(x,y) \} \ge 0], \;\forall y \in Y.
\end{alignat}
We have relaxed the constraints; thus $f^L$ is a lower bound on $f^*$.
It turns out that it is generically true that $f^* = f^L$;
that is, for ``most'' problem data in a certain class, we can expect the equality to hold.
More specifically, it is generically valid that the feasible set of \eqref{eq:gsip_rel} is the closure of the feasible set of \eqref{eq:gsip}.
See \cite{gunzelEA07}.

For simplicity, write
\[
\bar{h}(x,y) = \max_j\set{ h_j(x,y) } 
\]
so that
\[
 \bar{h}(x,y) \le 0 \iff h_j(x,y) \le 0, \forall j.
\]
Then we note that \eqref{eq:gsip_rel} may also be written as
\begin{alignat}{2} 
\notag
\inf_{x \in X} & f(x) \\
\st 				
\notag & \max\set{ g(x,y), \bar{h}(x,y) } \ge 0,\;\forall y \in Y,
\end{alignat}
which highlights that it is a ``standard'' SIP, with lower-level program (LLP)  
\begin{equation}
\tag{SIP LLP}
\label{eq:sip_llp}
\inf_y \set{ \max\set{ g(x,y), \bar{h}(x,y)} : y \in Y }
\end{equation}

The approach in \cite{mitsosEA15} to obtain a lower bound is, effectively, to solve the SIP \eqref{eq:gsip_rel} with something akin to the constraint-generation/discretization method of \cite{blankenshipEA76}.
The issue is that \eqref{eq:sip_llp} is not solved.
Solution of the correct lower-level program is fairly critical to the proof that the lower bounds generated by the discretization method converge to $f^L$.

The analysis in \cite{mitsosEA15} relies on the following assumptions.
\begin{assumption}~
\label{assm:all}
\begin{enumerate}\itemsep0pt \parskip0pt
\item \label{assm:all:cc}
The host sets $X$ and $Y$ are compact, and $f$, $g$, and $\bar{h}$ are continuous on them.
\item \label{assm:all:generic}
It holds that $f^* = f^L$.
\item \label{assm:all:slater}
Problem~\eqref{eq:gsip} is infeasible, or for a given $\epsilon > 0$ there exists $\delta > 0$ and $x^s \in X$ such that
\[
f(x^s) \le f^* + \epsilon
\quad\text{and} \quad
\left([g(x^s,y) \ge \delta] \lor [\bar{h}(x^s,y) \ge \delta]\right), \;\forall y \in Y.
\]
\end{enumerate}
\end{assumption}

Assumption~\ref{assm:all}.\ref{assm:all:cc} is a mild and standard assumption in global optimization.
Assumption~\ref{assm:all}.\ref{assm:all:generic} holds generically, as mentioned above.
Assumption~\ref{assm:all}.\ref{assm:all:slater} states that the problem is infeasible or else an $\epsilon$-optimal GSIP-Slater point exists.
This last assumption is not used in the analysis of the lower bounding method;
nevertheless the counterexamples in this work will satisfy all of these assumptions.

\section{Lower bounds}

The lower bounding method from \cite{mitsosEA15} can involve the solution of two different versions of lower-level programs, the original LLP of \eqref{eq:gsip} (see problem \eqref{eq:llp} below), and an auxiliary LLP.
We analyze the lower bounding method in two situations, beginning with the situation that only the GSIP LLP is solved.

\subsection{A sketch of the method and argument}

The setting of the method is the following.
The method is iterative and at iteration $k$, for a given finite subset $Y^{L,k} \subset Y$, a lower bound of $f^L$ is obtained from the finite program
\begin{alignat}{2} 
\label{eq:gsip_rel_lower}
f^{L,k} = 
\inf_{x}\; & f(x) \\
\st 				
\notag & x \in X, \\
\notag & [g(x,y) \ge 0] \lor [\bar{h}(x,y) \ge 0], \;\forall y \in Y^{L,k}
\end{alignat}
Let the minimizer be $x^k$.
We  can assume that the lower bounding problem \eqref{eq:gsip_rel_lower} is always feasible, otherwise we can conclude that \eqref{eq:gsip_rel} is infeasible (and the method would terminate in finite iterations).
In \cite{mitsosEA15}, the original GSIP LLP is solved.
We obtain
\begin{equation}
\tag{LLP}
\label{eq:llp}
y^k \in \arg\min_y\set{ g(x^k,y) : y \in Y, \bar{h}(x^k,y) \le 0 }.
\end{equation}
For the present analysis, assume that $\bar{h}(x^k,y^k) < 0$;
if this holds the lower bounding method does not require the solution of the auxiliary LLP (see the following section).
Furthermore, we may assume that $g(x^k,y^k) < 0$ 
(and, in particular, that \eqref{eq:llp} is feasible, or else we have found $x^k$ feasible in \eqref{eq:gsip_rel}, and in this case since $x^k$ is the global minimizer of the relaxation \eqref{eq:gsip_rel_lower}, it is also a global minimizer of \eqref{eq:gsip_rel}).
Then we set $Y^{L,k+1}  = Y^{L,k} \cup \set{y^k}$ and we iterate.

The claim is that the sequence of iterates $\seq[k\in \mbb{N}]{x^k}$ have an accumulation point $x^*$ which is feasible, and as a consequence, that  $f^{L,k}$ increases to $f^L$.
We will try (but ultimately fail) to prove this using the approach from \cite{blankenshipEA76} to better understand how a counterexample may be constructed.

Since the sequence $\seq[k]{(x^k,y^k)}$ is a subset of $X \times Y$ which is compact, a subsequence must converge to some point $(x^*,y^*)$
(abuse notation and denote this subsequence $\seq[k]{(x^k,y^k)}$).
First note the construction of $Y^{L,k}$ and $x^k$ implies that we have
\[
  [g(x^{\ell},y^k) \ge 0] \lor [\bar{h}(x^{\ell},y^k) \ge 0], \;\forall \ell > k,
\]
which follows from the feasibility of $x^{\ell}$ in \eqref{eq:gsip_rel_lower} and that $Y^{L,\ell} \supset \set{ y^k : k < \ell }$.
This is equivalent to 
$\max\set{g(x^{\ell},y^k),\bar{h}(x^{\ell},y^k)} \ge 0$, for all $\ell > k$.
Taking the limit over $\ell$ and  then $k$, we get 
\[
\max \set {g(x^*,y^*),\bar{h}(x^*,y^*)} \ge 0.
\]
Note that since $\max\set{g(x^k,y^k),\bar{h}(x^k,y^k)} < 0$ for all $k$, we can conclude in addition that
\[\max \set {g(x^*,y^*),\bar{h}(x^*,y^*)} = 0,\]
which implies that at least one of $g(x^*,y^*)$ or $\bar{h}(x^*,y^*)$ is zero.
This leads to two cases that will be analyzed shortly.

The next claim is that $x^*$ must be feasible in \eqref{eq:gsip_rel}.
If not, then there exists $y^{\dagger}\in Y$ with 
\[
\max \set{g(x^*,y^{\dagger}), \bar{h}(x^*,y^{\dagger})} < 0
\]
or that $g(x^*,y^{\dagger}) < 0$ and $\bar{h}(x^*,y^{\dagger}) < 0$.
By continuity, for all $k$ sufficiently large 
$g(x^k, y^{\dagger}) < 0$ and $\bar{h}(x^k, y^{\dagger}) < 0$.
By definition of $y^k$ (noting that $y^{\dagger}$ is feasible in \eqref{eq:llp} of which $y^k$ is the global minimizer),
$
g(x^k,y^k) \le g(x^k, y^{\dagger}) < 0,
$
so taking limits again
\[
g(x^*,y^*) \le g(x^*,y^{\dagger}) < 0.
\]
We can now analyze two cases.
First case: 
If $g(x^*,y^*) = 0$ (and $\bar{h}(x^*,y^*) \le 0$), we immediately have a contradiction.
Second case:
We must have $g(x^*,y^*) < 0$, and thus $\bar{h}(x^*,y^*) = 0$.
\emph{The issue is that we have not constructed $y^k$ as a global minimizer of $\bar{h}(x^k,\cdot)$, to arrive at a similar contradiction.}

\subsubsection*{Counterexample}

\begin{figure}
\begin{center}
\begin{tikzpicture}[xscale=1.5,yscale=1.5]
\draw[fill=gray] (-0.5,-1) -- (0.5,1) -- (1,1) -- (1,-1) -- cycle;
\draw[latex-latex] (0,-1.3) -- (0,1.3) node[above]{$y$}; 
\draw[latex-latex] (-1.3,0) -- (1.3,0) node[right]{$x$}; 
\draw (-1,-1) rectangle (1,1); 
\draw[dashed,domain=0:1] plot (\x, \x);
\end{tikzpicture}
\end{center}
\caption{Visualization of counterexample~\eqref{eq:counter_ex}.
The box represents $[-1,1] \times [-1,1]$.
The shaded grey area is the subset of $(x,y)$ such that $-2x + y \le 0$.
The dashed line represents the minimizers of the LLP for $x \ge 0$.}
\label{fig:cex1}
\end{figure}
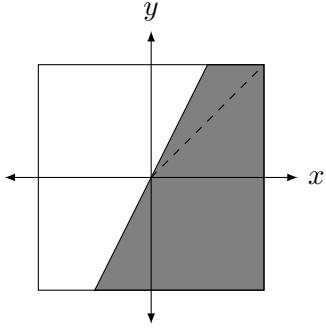

Consider 
\begin{alignat}{2} 
\tag{CEx 1}
\label{eq:counter_ex}
\inf_{x}\; & -x \\
\st 				
\notag & x \in [-1,1], \\
\notag & 0 \le (x-y)^2-10, \;\forall y \in [-1,1]: 
	-2x + y \le 0.
\end{alignat}
The behavior to note is this:
We are trying to maximize $x$;
For fixed $x$, the LLP feasible set consists of $y \le 2x$;
For $x < -\sfrac{1}{2}$, the LLP is infeasible;
For all $(x,y)$ in the domain $[-1,1]^2$, the LLP objective $g(x,y) = (x-y)^2-10$ is strictly less than zero;
The minimizer of the LLP at $x^k \ge 0$ is $y^k = x^k$;
The set of feasible $x$ is $[-1,-\sfrac{1}{2})$.
See Figure~\ref{fig:cex1}.

Clearly, then, the infimum is not attained at a feasible $x$, but the closure of the feasible set is indeed described by the SIP relaxation
\begin{alignat}{2} 
\notag
\inf_{x}\; & -x \\
\st 				
\notag & x \in [-1,1], \\
\notag & [0 \le (x-y)^2-10] \lor [-2x + y \ge 0], \;\forall y \in [-1,1].
\end{alignat}
We can ignore the always-false clause $0 \le (x-y)^2-10$.
Further, it is easy to see that the feasible set is
$\set{x \in [-1,1] : x \le (\sfrac{1}{2})y, \forall y \in[-1,1] } = [-1,-\sfrac{1}{2}]$
as hoped.
The infimum, consequently, is $\sfrac{1}{2}$.

Beginning with $Y^{L,1} = \emptyset$,  the minimizer of the lower bounding problem is $x^1 = 1$.
Solving the LLP, we get $y^1 = 1$ which we note satisfies the lower-level inequality strictly.
That is, $g(x^1,y^1) < 0$ and $\bar{h}(x^1,y^1) = -2x^1 + y^1 = -1 < 0$ as required by the lower bounding method in \cite{mitsosEA15}.

The next iteration, with $Y^{L,2} = \set{1}$, adds the constraint 
$-2x + 1 \ge 0$
to the lower bounding problem;
the feasible set is $[-1,\sfrac{1}{2}]$ so the minimizer is $x^2 = \sfrac{1}{2}$.
The minimizer of the LLP is $y^2 = \sfrac{1}{2}$;
again, it is in the interior of the feasible set and the optimal objective value is negative.

The third iteration, with $Y^{L,2} = \set{1, \sfrac{1}{2}}$, adds the constraint 
$-2x + \sfrac{1}{2} \ge 0$
to the lower bounding problem;
the feasible set is $[-1,\sfrac{1}{4}]$ so the minimizer is $x^3 = \sfrac{1}{4}$.
The minimizer of the LLP is $y^3 = \sfrac{1}{4}$;
again, it is in the interior of the feasible set and the optimal objective value is negative.

In general, we see that the iterates satisfy 
$x^k = y^k = \frac{1}{2^{k-1}}$.
Consequently, they converge to $0$.
In particular, the corresponding lower bounds $\seq[k]{ - x^k}$ converge to zero, which we note is strictly less than the infimum of $\sfrac{1}{2}$.

\subsection{Where it goes wrong}

The claim in  \cite{mitsosEA15} (which the counterexample shows is false) is that 
$g(x^*,y^*) \ge 0$
always holds.
The argument is that since $\bar{h}(x^k,y^k) < 0$ for all $k$, then there exists $K$ so that
$\bar{h}(x^{\ell},y^k) < 0$ for all $\ell \ge k \ge K$.
The conclusion is that for the disjunction in \eqref{eq:gsip_rel_lower} to be true, the clause $g(x^{\ell},y^k) \ge 0$  must hold for all sufficiently large $\ell,k$ such that $\ell > k$.
Taking the limit over $\ell$ and  then $k$, we get $g(x^*,y^*) \ge 0$.

However, the counterexample above demonstrates that 
$\bar{h}(x^k,y^k) = -2x^k + y^k = - \frac{1}{2^{k-1}} < 0$ for all $k$,
and yet for sufficiently large $\ell$ (in fact, for $\ell > k + 1$) we have
\[
\bar{h}(x^{\ell},y^k) = -2x^{\ell} + y^k =  - 2\frac{1}{2^{\ell-1}}  + \frac{1}{2^{k-1}} > 0.
\]
The claim that there exists $K$ so that $\bar{h}(x^{\ell},y^k) < 0$ for all $\ell > k > K$ is false.
Consequently, we cannot conclude that $g(x^*,y^*) \ge 0$ always holds.

\subsection{Analysis with auxiliary LLP}

The method in \cite{mitsosEA15} introduces the ``auxiliary LLP'' which is attempting to get closer to \eqref{eq:sip_llp}:
\begin{equation}
\tag{AUX LLP}
\label{eq:aux_llp}
\tilde{y}^k \in \arg \min_y \set{ \bar{h}(x^k,y) : y \in Y, g(x^k,y) \le \alpha g(x^k,y^k)}
\end{equation}
with $\alpha \in (0,1)$.
That is, we minimize the constraints subject to being approximately optimal
(recall that $y^k$ is the global minimizer of \eqref{eq:llp}).
The method as stated does not require that it is always solved (and the counterexample takes advantage of this).
Unfortunately, even if the auxiliary LLP is always solved, we do not get correct behaviour of the method.

We analyze the same method as before but instead populate $Y^{L,k}$ with the minimizers of the auxiliary LLP:
\[
Y^{L,k+1}  = Y^{L,k} \cup \set{\tilde{y}^k}.
\]
The analysis of the claim that $f^{L,k}$ increases to $f^L$ proceeds similarly to before;
we can assume $g(x^k,y^k) < 0$ (or else we have found a feasible and thus global optimal point of \eqref{eq:gsip_rel}),
and thus $g(x^k,\tilde{y}^k)< 0$.
We can also assume that \eqref{eq:llp} is feasible, and so
$\bar{h}(x^k,\tilde{y}^k) \le \bar{h}(x^k,y^k) \le 0$.
Then, as before, we can conclude that there is a (sub)sequence of iterates 
$\seq[k]{(x^k,\tilde{y}^k)}$ 
converging to $(x^*,y^*)$ with
$\max \set {g(x^*,y^*),\bar{h}(x^*,y^*)} \ge 0$
and again, since 
$g(x^k,\tilde{y}^k) < 0$ and
$\bar{h}(x^k,\tilde{y}^k) \le 0$
this implies 
$\max \set {g(x^*,y^*),\bar{h}(x^*,y^*)} = 0$
and so at least one of $g(x^*,y^*)$ or $\bar{h}(x^*,y^*)$ is zero.

Again, if $x^*$ is not feasible in \eqref{eq:gsip_rel}, then there exists $y^{\dagger}\in Y$ with  $g(x^*,y^{\dagger}) < 0$ and $\bar{h}(x^*,y^{\dagger}) < 0$.
If $g(x^*,y^*)$ is zero, then we may be able to derive a contradiction.
However, in the other case that $g(x^*,y^*) < 0$, and $\bar{h}(x^*,y^*) = 0$, we cannot.
This is because it is possible that 
\[
\alpha g(x^k,y^k) < g(x^k,y^{\dagger}) < 0
\]
for all sufficiently large $k$.
This means that $y^{\dagger}$ is not feasible in the auxiliary LLP, and so this allows the possibility that
\[
\bar{h}(x^k,y^{\dagger}) \le \bar{h}(x^k,\tilde{y}^k)
\]
with
$\bar{h}(x^k,y^{\dagger}) \to \bar{h}(x^*,y^{\dagger}) < 0$
and
$\bar{h}(x^k,\tilde{y}^k) \to \bar{h}(x^*,y^{*}) = 0$.
This is demonstrated with the following counterexample.

\subsubsection*{Counterexample 2}

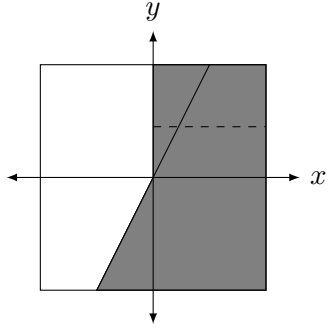
\begin{figure}
\begin{center}
\begin{tikzpicture}[xscale=1.5,yscale=1.5]
\draw[fill=gray] (-0.5,-1) -- (0,0) -- (0,1) -- (1,1) -- (1,-1) -- cycle;
\draw (-0.5,-1) -- (0.5,1); 
\draw[latex-latex] (0,-1.3) -- (0,1.3) node[above]{$y$}; 
\draw[latex-latex] (-1.3,0) -- (1.3,0) node[right]{$x$}; 
\draw (-1,-1) rectangle (1,1); 
\draw[dashed] (0,0.45) -- (1,0.45);
\end{tikzpicture}
\end{center}
\caption{Visualization of counterexample~\eqref{eq:counter_ex2}.
The box represents $[-1,1] \times [-1,1]$.
The shaded grey area is the subset of $(x,y)$ such that $\min\set{-2x + y, -x} \le 0$.
The dashed line represents $y = 0.45$ and helps visualize the feasible set of the auxiliary LLP.}
\label{fig:cex2}
\end{figure}

We have a similar setting as before.
Consider 
\begin{alignat}{2} 
\tag{CEx 2}
\label{eq:counter_ex2}
\inf_{x}\; & -x \\
\st 				
\notag & x \in [-1,1], \\
\notag & 0 \le -y-10, \;\forall y \in [-1,1]: 
	\min\set{-2x + y, -x} \le 0.
\end{alignat}
The behavior to note is this:
We are trying to maximize $x$;
For fixed $x$, the LLP feasible set consists of
$y \le 2x$ 
OR 
all $y \in [-1,1]$ if $x \ge 0$;
For $x < -\sfrac{1}{2}$, the LLP is infeasible;
For all $(x,y)$ in the domain $[-1,1]^2$, the LLP objective $g(x,y) = -y-10$ is strictly less than zero;
The set of feasible $x$ is $[-1,-\sfrac{1}{2})$.
Again, the infimum is $\sfrac{1}{2}$.
See Figure~\ref{fig:cex2}.

The SIP relaxation is
\begin{alignat}{2} 
\notag
\inf_{x}\; & -x \\
\st 				
\notag & x \in [-1,1], \\
\notag & [0 \le -y-10] \lor [\min\set{-2x + y, -x} \ge 0], \;\forall y \in [-1,1].
\end{alignat}
We can ignore the always-false clause $0 \le -y-10$.
Further, it is easy to see that the feasible set is
\[
\set{x \in [-1,1] : [x \le 0] \land [x \le (\sfrac{1}{2})y, \forall y \in[-1,1]] }
\]
which equals $[-1,-\sfrac{1}{2}]$.
The infimum is $\sfrac{1}{2}$.

Beginning with $Y^{L,1} = \emptyset$, the minimizer of the lower bounding problem is $x^1 = 1$.
The minimizer of the LLP is $y^1 = 1$ with objective value $-11$.
If we choose $\alpha = 0.95$, then the feasible set of the auxiliary LLP is 
$
\set{ y \in [-1,1] : -y-10 \le -10.45}
$
or $y \in [0.45,1]$.
The minimizer of the auxiliary LLP is $\tilde{y}^1 = 0.45$;
the optimal objective value is 
$\min\set{-2x^1 + \tilde{y}^1,  -x^1} = \min\set{-2 + 0.45,  -1} = -1.55$.

In the next iteration, with $Y^{L,2} = \set{0.45}$, the lower bounding problem has feasible set
\[
\set{x \in [-1,1] : [x \le 0] \land [x \le 0.225] } = [-1,0]
\]
and so yields $x^2 = 0$.
The minimizer of the LLP is again $y^2 = 1$ with objective value $-11$.
The feasible set of the auxiliary LLP is again $y \in [0.45,1]$.
The objective function of the auxiliary LLP is $\min\set{-2x^2 + y, -x^2} = 0$ for all feasible $y$, and so all feasible $y$ are optimal.

However, in the third iteration, no matter what the value of $\tilde{y}^2$ is, the lower bounding solution is again $x^3 = 0$.
This is because the feasible set is
\[
\set{x \in [-1,1] : [x \le 0] \land [x \le 0.225] \land [x \le \tilde{y}^2/2] } = [-1,0]
\]
for any possible value of $\tilde{y}^2 \in [0.45,1]$.
The LLP and auxiliary LLP are the same, and again $\tilde{y}^3$ does nothing to change the feasible set of the lower bounding problem.

The sequence $\seq[k]{x^k}$ trivially converges to $0$, and the lower bounds also converge to $0$, which is strictly less than the infimum of $\sfrac{1}{2}$.

\section{Remarks}

The lower bounds for Example~\eqref{eq:counter_ex2} may converge with a different value of $\alpha$, or a non-empty initialization of $Y^{L,1}$,  but neither of these guarantee the general behaviour of the lower bounding method claimed in \cite{mitsosEA15}.
It seems that a convergent lower bounding method relies on the solution of \eqref{eq:sip_llp} and populating $Y^{L,k}$ with its minimizers.
This characterizes recent approaches like in \cite{djelassiEA19}.

\bibliographystyle{plain}

\end{document}